\documentclass[12pt]{article}
\usepackage{amssymb,amsmath,amsfonts,mathrsfs,upgreek} 
\usepackage{mathabx}

\usepackage{graphicx}

\usepackage[colorlinks=true]{hyperref}
\usepackage{comment}

\renewcommand\comment[1]{{\iffalse #1 \fi}}

\newtheorem{theorem}{Theorem}

\newtheorem{assumption}{Assumption}[section]

\newtheorem{remark}{Remark}
\newtheorem{notation}{Notation}[section]
\newtheorem{definition}{Definition}[section]

\newcommand{\io}{{\infty}}

\newcommand{\real}{ {\mathbb R}   }
\newcommand{\torus}{ {\mathbb T}    }
\newcommand{\integer}{ {\mathbb Z}   }

\newcommand{\complex}{ {\mathbb C}   }
\newcommand{\bB}{ {\mathbb B}   }

\newcommand{\cP}{{\mathbb P}^n_{\! s}}

\newcommand\beq[1]{ \begin{equation}\label{#1} }
\newcommand{\eeq}{ \end{equation} }
\newcommand{\beqno}{ \[ }
\newcommand{\eeqno}{ \] }
\newcommand\beqa[1]{ \begin{eqnarray} \label{#1}}
\newcommand{\eeqa}{ \end{eqnarray} }
\newcommand{\beqano}{ \begin{eqnarray*} }
\newcommand{\eeqano}{ \end{eqnarray*} }

\newcommand\dfn[1]{ \begin{definition}\label{#1} }
\newcommand\edfn{ \end{definition} }
\newcommand\ass[1]{ \begin{assumption}\label{#1} }
\newcommand\eass{ \end{assumption} }

\newcommand\notat[1]{ \begin{notation} \label{#1} 
 }
\newcommand\enotat{\end{notation}}

\newcommand\rem{\begin{remark} 
\rm 
}
\newcommand\erem{\end{remark} %\vglue0.7truecm\noindent
}

\newcommand{\nl}{{\smallskip\noindent}}
\newcommand{\giu}{{\medskip\noindent}}

\newcommand{\e}{\varepsilon}

\renewcommand{\a }{\alpha }

\newcommand{\ii }{{\rm i} }

\newcommand{\g }{\gamma}

\renewcommand{\t }{\tau }
\renewcommand{\o }{\omega }
\renewcommand{\O }{\Omega }

\newcommand{\Z}{\mathbb{Z}}

\renewcommand\varpi{\upkappa}

\newcommand{\Nf}{\mathtt N}

\def\R{\mathbb R}
\def\T{\mathbb T}

\def\dst{\displaystyle}
\def\bks{\, \backslash\, }
\def\meas{{\rm\, meas\, }}
\newcommand\ham{\mathtt{H}}

\newcommand\sa{\theta} %\newcommand\sa{\theta}
\newcommand{\fproj}{\pi_{\!{}_{\integer k}}}

\usepackage{appendix}

\newcommand\gen{{\cal G}^n}

\newcommand\noruno[1]{  |#1|_{{}_1} }

\renewcommand\ln{\log}
\newcommand\st{{\rm \ s.t.\ }}
\newcommand\hol{{\bB}}

\renewcommand\subset{\subseteq}
\newcommand\DD{\mathcal D}

\renewcommand\O{\DD}

\makeatletter
\newcommand{\pushright}[1]{\ifmeasuring@#1\else\omit\hfill$\displaystyle#1$\fi\ignorespaces}
\newcommand{\pushleft}[1]{\ifmeasuring@#1\else\omit$\displaystyle#1$\hfill\fi\ignorespaces}
\makeatother

\title{\bf 
Quasi--periodic motions in generic nearly--integrable mechanical systems
}

\begin{document}

\author{ 
%{\bf Scientific chapter:} {\sl Mathematical analysis.} \\
\footnotesize L. Biasco  \& L. Chierchia
\\ \footnotesize Dipartimento di Matematica e Fisica, Universit\`a degli Studi Roma Tre
\\ \footnotesize Largo San L. Murialdo 1 - 00146 Roma, Italy
\\ {\footnotesize biasco@mat.uniroma3.it, luigi.chierchia@uniroma3.it}
\\ 
}

\maketitle

%\date{}

\begin{abstract}\noindent
In this note we present and briefly discuss results, which include as a particular case the theorem announced  in \cite{BClin}, concerning the typical behaviour of nearly--integrable mechanical systems with  generic analytic potentials.
\end{abstract}

%\section{Introduction}

\nl
In 2015, encouraged by our mentor, colleague and friend Antonio Ambrosetti, we published in {\sl Atti della Accademia Nazionale dei Lincei}
 an announcement \cite{BClin} concerning `typical'  trajectories of nearly--integrable Hamiltonian systems. In particular, we stated a theorem (\cite[p. 426]{BClin}), which,  can be roughly rephrased as follows: 

\nl
{\sl In bounded regions of  phase space, except for a set of measure $\e|\log\e|^\g$, trajectories of  nearly--integrable mechanical systems on $\real^n\times \torus^n$ with generic real--analytic potentials of size $\e\ll 1$ are quasi--periodic and span $n$--tori invariant for such systems.} 

\nl
This theorem 
is in agreement (up to the logarithmic correction) with a conjecture formulated by Arnold, Kozlov and Neishtadt in the Springer 
Encyclopaedia of Mathematical Sciences \cite[Chapter 6, p. 285]{AKN}. 

\nl
A complete proof of the above result turned out to be much longer and more delicate than we thought and it has been completed only recently in  \cite{BCaa} and \cite{BCtheend},  which in turn exploit intermediate results published in \cite{BCKAM} and  \cite{BCnonlin}.

\nl
The purpose of this short note is to communicate  the precise results of \cite{BCtheend}, which,  as a particular case, yield the above theorem.

%\section{Results}

\nl

\nl
In order to state the main results in \cite{BCtheend} we need to recall a few notions and give some definitions.

\begin{itemize}

\item[\rm (a)] {\sl Hamiltonian systems on $\real^n\times \T^n$}
\\
Given 
a region $B\subset \real^n$, the `phase space' ${\mathcal M}:=B\times \torus^n$ (where $\T^n:=\real^n/(2\pi \Z^n)$) and a real analytic `Hamiltonian function' $H:{\cal M}\to\real$,
we denote by $z\in{\cal M}\to  \Phi^t_H(z)\in{\cal M}$  the {\sl Hamiltonian flow generated by $H$}, namely, the solution of the standard Hamilton equations
$$
\left\{
\begin{array}{l}
\dot y = -H_x(y,x)\\
\dot x= H_y(y,x)
\end{array}\right.\,,\qquad (y,x)|_{t=0}=z\,,
$$
where, as usual, `dot' denotes derivative with respect to `time' $t\in\real$, and $H_y$, $H_x$ the gradient with respect to $y$, $x$. \\
A {\sl mechanical system} on $\real^n\times \torus^n$ is a Hamiltonian system with Hamiltonian 
$${\rm H}(y,x)=\frac12 |y|^2+f(x)\,,\qquad \big({\rm where}\ |y|^2:=y\cdot y:=\sum_j|y_j|^2\big)\,,
$$
whose evolution equations are equivalent to the Newton equation $\ddot x=-f_x(x)$; $f$ is called the {\sl potential} of the system; `nearly--integrable' means that 
the potential is of the form $\e f$ with $\e$ a small real parameter.

\item[\rm (b)] {\sl Diophantine vectors}\\
A vector $\o\in\real^n$ is called  {\sl Diophantine} if there exist $\a>0$ and $\t\ge n-1$ such that $\dst |\o\cdot k| \ge \a/|k|_{{}_1}^\t$, for any non vanishing integer vector $k\in\integer^n$, where $\noruno{k}:=\sum |k_j|$.

\item[\rm (c)] {\sl Maximal KAM tori}\\
A  set ${\mathcal T}\subset {\cal M}$ is a {\sl maximal KAM torus} for a Hamiltonian function $H$ if there exist a real analytic embedding $\phi:\torus^n\to {\cal M}$ and a Diophantine frequency vector $\o\in\real^n$ such that ${\mathcal T}=\phi(\torus^n)$ and for each $z\in{\mathcal T}$, $\Phi^t_H(z)=\phi(x+ \o t)$, where
$x=\phi^{-1}(z)$. For general information on KAM (Kolmogorov, Arnold, Moser) Theory, see \cite{AKN}, and references therein.

\item[\rm (d)] {\sl Generators of 1d maximal lattices}\\
Let $ \integer^n_\varstar$ be the set of integer vectors $k\neq 0$ in $\integer^n$ such that the  first non--null  component is positive:
\beqno
 \integer^n_\varstar:=
 \big\{ k\in\integer^n:\ k\neq 0\ {\rm and} \ k_j>0\ {\rm where}\ j=\min\{i: k_i\neq 0\}\big\}\,.
 \eeqno
$\gen$  denotes the set of {\sl generators of 1d maximal lattices} in $\integer^n$, namely, the set of  vectors $k\in  \integer^n_\varstar$ such that the greater common divisor (gcd)  of their components is 1:
\beqno
\gen:=\{k\in \integer^n_\varstar:\ {\rm gcd} (k_1,\ldots,k_n)=1\}\,.
\eeqno

\item[\rm (e)] {\sl Resonances}\\
A resonance ${\cal R}_k$
with respect to the free Hamiltonian $\frac12 |y|^2$ is the set $\{y\in \real^n: y\cdot k=0\}$, where $k\in\gen$. 
We call
${\cal R}_{k,\ell}$  a {\sl double resonance} if  ${\cal R}_{k,\ell}={\cal R}_k\cap {\cal R_\ell}$ with $k$ and $\ell$ in $\gen$ linearly independent; the order of a double resonance is given by $\max\{\noruno{k},\noruno{\ell}\}$.

\item[\rm (f)] {\sl 1d Fourier projectors}\\
Given $k\in \integer^n\bks\{0\}$ and a periodic analytic function $f:\torus^n\to\complex$, we denote by $\fproj f$ the (analytic) periodic function of  
{\sl one variable} $\sa\in\torus$ given by
$$
\sa\in\torus\mapsto \fproj f (\sa):=\sum_{j\in\integer} f_{jk} e^{\ii j\sa}\,.
$$
Note that $\dst f(x)= \sum_{k\in\gen}  \fproj f  (k\cdot x)$.

\item[\rm (g)] {\sl Morse functions with distinct critical values}\\
A  function $\sa\to F(\sa)$ is a Morse function if its critical points are non--degenerate, i.e., $F'(\sa_0)=0\implies F''(\sa_0)\neq 0$; `distinct critical values' means that if 
$\sa_1\neq\sa_2$  are distinct critical points, then $F(\sa_1)\neq F (\sa_2)$.

\item[\rm (h)] {\sl A Banach space of real analytic functions}\\
Let $s>0$.  We denote by $\hol_s^n$   the Banach space of
real analytic periodic functions on $\torus^n$ 
having zero average:
\beqno
\hol_s^n:=\Big\{f=\sum_{k\in\Z^n \atop k\neq 0} f_k e^{\ii k\cdot x} \st \bar f_{k}=f_{-k}\ {\rm and }\ \ \|f\|_s<\io\Big\}\,,
\eeqno
where $\dst \|f\|_s:=\sup_{k\in \integer^n}  |f_k| e^{|k|_{{}_1}s}$.

\item[\rm (i${}^\star$)] {\sl The generic set $\cP$ of potentials}\\
We  denote by $\cP$ the subset of the unit ball of $\hol_s^n$ given by the set of functions $f\in \hol_s^n$ such that the following two conditions 
hold:
\begin{flalign}
\label{P1}
&\phantom{.} \varliminf_{\noruno{k}\to+\io\atop k\in\gen} |f_k| e^{\noruno{k}s} \noruno{k}^n>0\,,
\nonumber\\ 
\nonumber
&\phantom{a.} 
\forall \ k\in\gen, \fproj f\ {\rm is \, a\, Morse\, function\, with \, distinct \, critical\, values}.
\end{flalign}

\end{itemize}
We remark that all the above definitions are standard, except for the last one, which describe the class of potentials for which our results hold. \\
   $\cP$ is a typical set
in many ways: it contains an open and dense set (in the topology of~$\hol^n_s$), it has full measure with respect to standard probability measures on the unit ball of $\hol^n_s$, and is a prevalent set (`prevalence' is a measure--theoretic notion for subsets of infinite--dimensional spaces that is analogous to `full Lebesgue measure' in Euclidean spaces; compare \cite{HK}). For a detailed  discussion of the properties of $\cP$, see, 
\cite[\S 3]{BCnonlin} and Appendix~A.2 in \cite{BCtheend}. 
%:
We also remark that the definition given here simplifies and extends  former definitions given in \cite{BClin} and \cite{BCnonlin}.

\nl
We can now state the main results in \cite{BCtheend}.

\begin{theorem}\label{prometeo1} Let $n\ge 2$, $s>0$, $0<\e<1$, $f\in \cP$, $B$ an open ball in $\R^n$ and  $\ham(y,x;\e):=\frac12 |y|^2 +\e f(x)$. 
Then, there exists a constant $c>1$ such that
   all  points in $B\times \T^n$
 lie on a maximal KAM torus for $\ham$, except for 
   a subset of measure bounded by
 $
 \, c\, \e |\ln \e|^{\g}$ with $\g:=11 n +4$. 
\end{theorem}

\begin{theorem}\label{prometeo2}
Fix  $0<a<1$. For any $\e>0$, 
there exists an
 open neighbourhood $\O^2\subseteq B$
 of double resonances of order smaller than $1/\e^{b}$,
 with $b:=\frac{1-a}{\g}$, which satisfies  
 $
\meas
\big(\O^2\times\T^n\big)\leq {\rm c_o} \e^a$,
for a suitable constant $\rm c_o$ (depending only on $n$),
such that the following holds. Under the assumptions of Theorem~\ref{prometeo1}, 
there exists  a positive constant  $\hat c$
 (independent of $a$) such that
 all points in $(B\setminus\O^2)\times \T^n$,
 lie on a maximal KAM torus for $\ham$,
except for an exponentially
small subset of measure bounded by 
$
e^{-\hat c/\e^b}
$.
\end{theorem}

\begin{theorem}\label{prometeo3} Let the assumptions of Theorem~\ref{prometeo1} hold and let $n=2$.
There exists a  constant $\bar c>0$
  such that, for every $0<a<1$, 
   all  points in $\{y\in B: |y|>\e^{a/2}\}\times\T^2$
 lie on a maximal KAM torus for $\ham$, except for
 an exponentially
small subset of measure bounded by 
$
e^{-\bar c/\e^b}
$,
with $b=(1-a)/24$.
\end{theorem}
Let us a make a few observations.

\giu
Theorem 1 -- which extends the result in \cite{BClin} --  may be viewed as the `ultimate frontier of KAM Theory', in the sense that, as remarked by Arnold et al., near double resonances, 
there are regions of order $\e$ where the dynamics of $\ham$ is equivalent to the dynamics of the parameter--free Hamiltonian
$\frac12 |y|^2 +f(x)$ and, therefore, 
{\sl it is natural to expect that in a generic system with three or more degrees of freedom the measure of the ``non-torus'' set has order $\e$}
 (\cite[Remark 6.18, p. 285]{AKN}). 
Theorem~1 provides an upper bound on the measure of the non--torus set
in agreement (up to the logarithmic correction $|\log\e|^\g$) 
with this expectation.
On the other hand, rigorous  lower bounds on such a measure
appear to be extremely hard to be  proven in the analytic case; for partial results in the Gevrey case, see~\cite{LMS}.

\giu
The KAM tori constructed in Theorem~1 are {\sl not} uniformly distributed in phase space. Indeed, if one stays away from a finite number of double resonances, the density is exponentially small: this is the main content of Theorem~2.

\giu
Theorem 3 is a consequence  of Theorem 2, since in dimension 2, the only double resonance in a mechanical system is the origin. Theorem~3 is in agreement with the conjecture formulated by Arnold et al. in \cite[Remark 6.17, p. 285]{AKN}.\\ A related (weaker) result was announced in  \cite{BC-DCDS}.

\giu
We recall that classical KAM Theory yields only {\sl primary tori} (which are graphs over $\T^n$) $\sqrt\e$--away from  resonances, while the new tori constructed in the above theorems fill, with an exponential density, a neighbourhood of  (simple) resonances  far from double resonances. Furthermore,  the new KAM tori include, besides primary tori, also {\sl secondary} tori, which exhibit different topologies and, in particular, are not  graphs over $\T^n$. \\
Secondary tori close to  resonances have been also investigated in \cite{MNT}.

\giu
To prove the above results  it is essential to study regions  close   to   resonances ${\cal R}_k$ for $\noruno{k}\to \infty$   as $\e\to0$. Away from  doubles resonances the `secular' dynamics in  ${\cal R}_k$ is, up to exponentially long times,  ruled by the integrable Hamiltonian
$$\ham_k:= \frac12 |y|^2 +\e (\fproj f)(x\cdot k)\,.$$ 
Therefore, it  should not surprise that one needs non--degeneracy assumptions in (i$^\star$) above; in particular, the first condition implies that for $\noruno{k}\ge \Nf$ large enough, but {\sl   independent of $\e$}, the secular potential $\fproj f$ is essentially a rescaled and shifted cosine (as fully discussed in \cite{BCnonlin}). Notice that for low modes ($\noruno{k}< \Nf$) the secular potential $\fproj f$ is a generic periodic function. In particular,  the phase portrait of $\ham_k$ is quite arbitrary and may have an arbitrary number of equilibria and separatrices. 
The main point here  is  to prove 
the persistence of {\sl all}  
 integrable tori of $\ham_k$ up to an exponentially small set (away from double resonances).

\giu 
We finally remark that one of the main issues in the proof of measure estimates
 is to show that the integrable secular Hamiltonian $\ham_k$ above, {\sl in its action  variables}, is Kolmogorov non--degenerate,
namely the action--to--frequency map is invertible. While  
 $\ham_k|_{{}_{\e=0}}=\frac12 |y|^2$ is obviously non--degenerate, 
it is a fact that this is not always true for $\ham_k$ (in its action variables) when  $\e>0$.
Indeed, 
 this is a {\sl singular perturbation problem}, as suggested by the fact that  
the level sets of   $\ham_k$ have different topologies,
due to the presence of secondary tori for $\e>0$.

%%%%%%%%%%%%%%%%%%%%%%%%%%%%%%%%%%%%%%%

\end{document}